\documentclass{article}
\usepackage{amssymb,amsmath}
\usepackage{anysize}
\usepackage{psfrag}
\usepackage{url}
\newtheorem{thm}{Theorem}[section]
\newtheorem{def.}{Definition}[section]

\newtheorem{prop}{Proposition}[section]
\newtheorem{cor}{Corollary}[section]

\numberwithin{table}{section}

\begin{document}
\title{Equivalence Classes of Colorings}
\author{Jun Ge\\
        School of Mathematical Sciences\\
        Xiamen University\\
        Xiamen, Fujian 361005\\
        P. R. China\\
        \texttt{mathsgejun@163.com}\\
        and\\
        Slavik Jablan\\
        The Mathematical Institute\\
        Knez Mihailova 36\\
        P.O. Box 367, 11001, Belgrade\\
        Serbia\\
        \texttt{sjablan@gmail.com}\\
        and\\
        Louis H. Kauffman\\
        Department of Mathematics, Statistics and Computer Science\\
        University of Illinois at Chicago\\
        851 S. Morgan St., Chicago IL 60607-7045\\
        USA\\
        \texttt{kauffman@uic.edu}\\
        and\\
        Pedro Lopes\\
        Center for Mathematical Analysis, Geometry and Dynamical Systems\\
        Department of Mathematics\\
        Instituto Superior T\'ecnico, Universidade de Lisboa\\
        Av. Rovisco Pais\\
        1049-001 Lisbon\\
        Portugal\\
        \texttt{pelopes@math.ist.utl.pt}\\
}
\date{September 17, 2013}
\maketitle

\bigbreak

\begin{abstract}
For any link and for any modulus $m$ we introduce an equivalence relation on the set of non-trivial $m$-colorings of the link (an $m$-coloring has values in $\mathbf{Z}/m\mathbf{Z}$). Given a diagram of the link, the equivalence class of a non-trivial $m$-coloring is formed by each assignment of colors to the arcs of the diagram that is obtained from the former coloring by a permutation of the colors in the arcs which preserves the coloring condition at each crossing. This requirement implies  topological invariance of the equivalence classes. We show that for a prime modulus the number of equivalence classes depends on the modulus and on the rank of the coloring matrix (with respect to this modulus).
\end{abstract}

\bigbreak


\bigbreak


\bigbreak

\section{Introduction} \label{sect:intro}

\noindent

\begin{figure}[!ht]
	\psfrag{ai}{\huge$c_i$}
	\psfrag{ai+1}{\huge$c_{i+1}$}
	\psfrag{aji}{\huge$c_{j_i}$}
	\psfrag{eq}{\huge$c_i + c_{i+1} - 2c_{j_i} = 0$}
	\centerline{\scalebox{.5}{\includegraphics{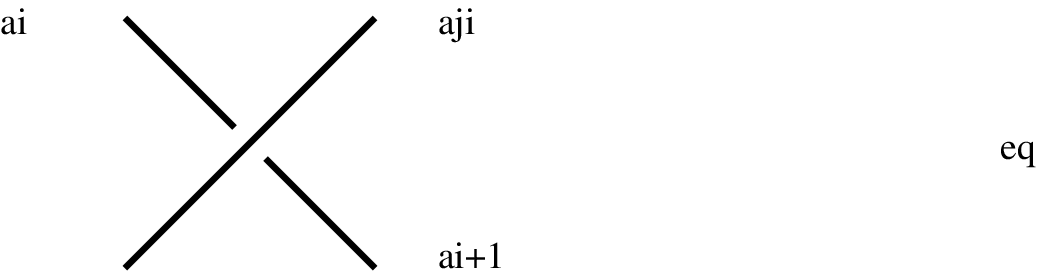}}}
	\caption{Arcs at a crossing and the equation read off it. The coloring system of equations is formed by each of these equations, one per crossing of the diagram under study.}\label{fig:xtop}
\end{figure}

Given a diagram $D$ of a link and a modulus $m>1$, a (Fox) coloring (\cite{CFox, lhKauffman}) is an assignment of integers modulo $m$ to the arcs of $D$ such that at each crossing twice the color assigned to the over-arc equals the sum of the colors assigned to the under-arcs, modulo $m$ (see Figure \ref{fig:xtop}).
For each diagram and for each modulus $m>1$ there is always at least one solution to this problem namely by assigning the same color (i.e., integer modulo $m$) to each and every arc of the diagram; thus there are exactly $m$ such solutions modulo $m$. These are the trivial solutions modulo $m$ i.e., the so-called trivial $m$-colorings of the diagram. The non-trivial $m$-colorings are the solutions, modulo $m$, which involve at least two distinct colors.

\bigbreak

{\bf Remark}. We remark that it is well known that this system of equations is also a system of relations for the first homology group of the 2-fold branched covering along the link (\cite{Przytycki}, Theorem 3.3). In fact, the fundamental group of the $2$-fold branched covering along a link is presented by labeling the arcs of the unoriented link diagram and having relations of the form $c=ba^{-1}b$ read off at each crossing when $b$ is the label of the over-crossing line. It then follows that $H_1(M_{2}(L))\oplus \mathbf{Z}$ (the first homology group of the $2$-fold branched covering along the link $L$) has presentation with $C=B-A+B=2B-A$, where $A, B, C$ are the corresponding elements in the abelianization of the fundamental group (\cite{Przytycki, Wada, Winker}). Should one set the color of one of the arcs equal to $0$ then there would be a bijective correspondence between this set of colorings and $H_1(M_{2}(L))$. It is interesting to remark that the fundamental group of the $2$-fold branched covering along the link is itself a non-abelian generalization of the Fox coloring. While we do not use this aspect of the topology here, we are aware of it and it may be of use in later work. For more background on this material see \cite{Fox, lhKauffmanOnKnots, Reidemeister, Seifert}.

\bigbreak

If a diagram endowed with an $m$-coloring undergoes a Reidemeister move, there is a unique reassignment of colors to the arcs involved in the move such that the new assignment is an $m$-coloring of the resulting diagram. Since these reassignments are reversible there is a bijection between the $m$-colorings before and after the performance of a finite number of Reidemeister moves. Furthermore, these reassignments preserve trivial $m$-colorings and thus they preserve also non-trivial $m$-colorings.

Therefore the number of $m$-colorings is a link invariant; the fact that a diagram of a link admits or not non-trivial $m$-colorings is an invariant of that link. It is known that there are links which do not admit non-trivial colorings over a given modulus. For example, the trefoil only admits non-trivial colorings over  moduli divisible by $3$.

\begin{figure}[!ht]
	\psfrag{0}{\huge$0$}
	\psfrag{1}{\huge$1$}
	\psfrag{2}{\huge$2$}
	\psfrag{3}{\huge$3$}
	\psfrag{4}{\huge$4$}
    \psfrag{x}{\huge}
\centerline{\scalebox{.5}{\includegraphics{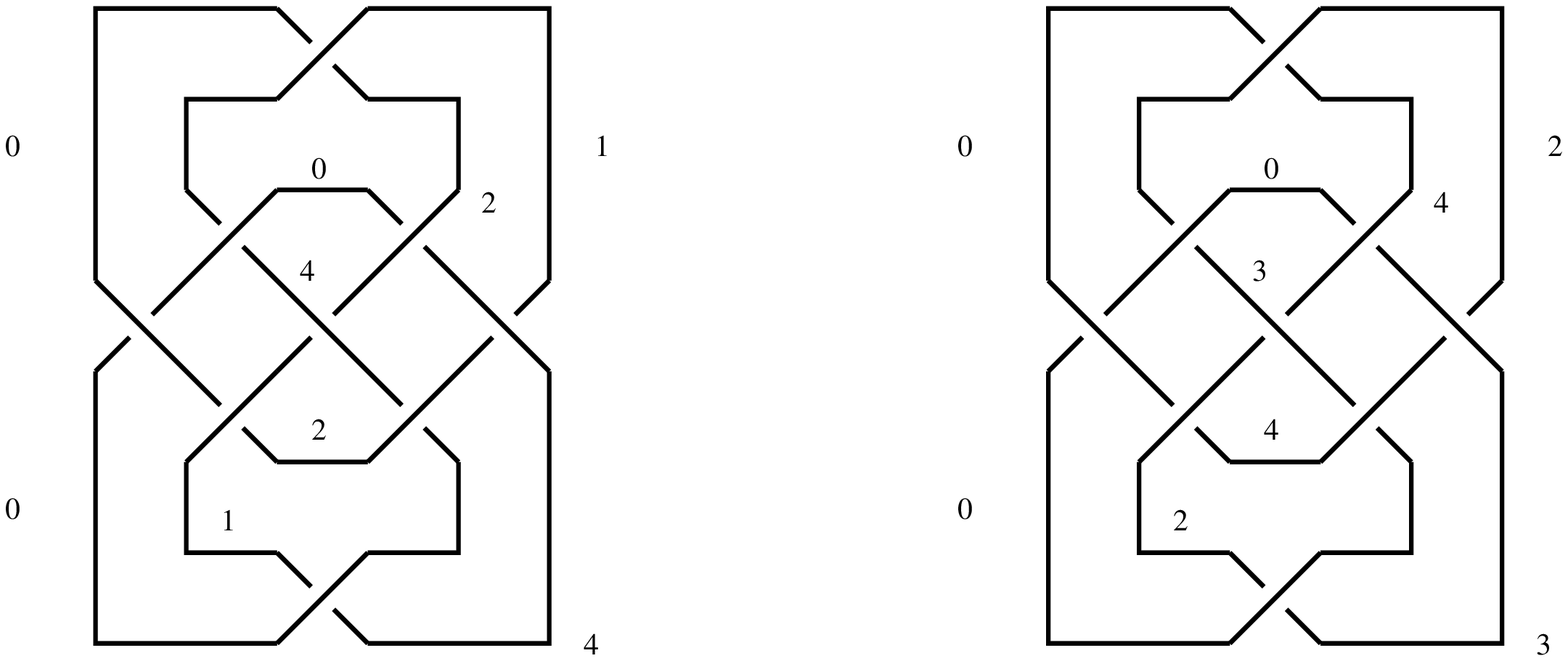}}}
	\caption{Two identical diagrams of $9_{40}$ but endowed with distinct non-trivial $5$-colorings. However, the colors of the one on the right are obtained from those of the one on the left by applying the permutation $(0)(1\, 2\, 4\, 3)$. These two $5$-colorings are related.}\label{fig:9_40tri}
\end{figure}

\begin{figure}[!ht]
	\psfrag{0}{\huge$0$}
	\psfrag{1}{\huge$1$}
	\psfrag{2}{\huge$2$}
	\psfrag{3}{\huge$3$}
	\psfrag{4}{\huge$4$}
    \psfrag{x}{\huge}
\centerline{\scalebox{.5}{\includegraphics{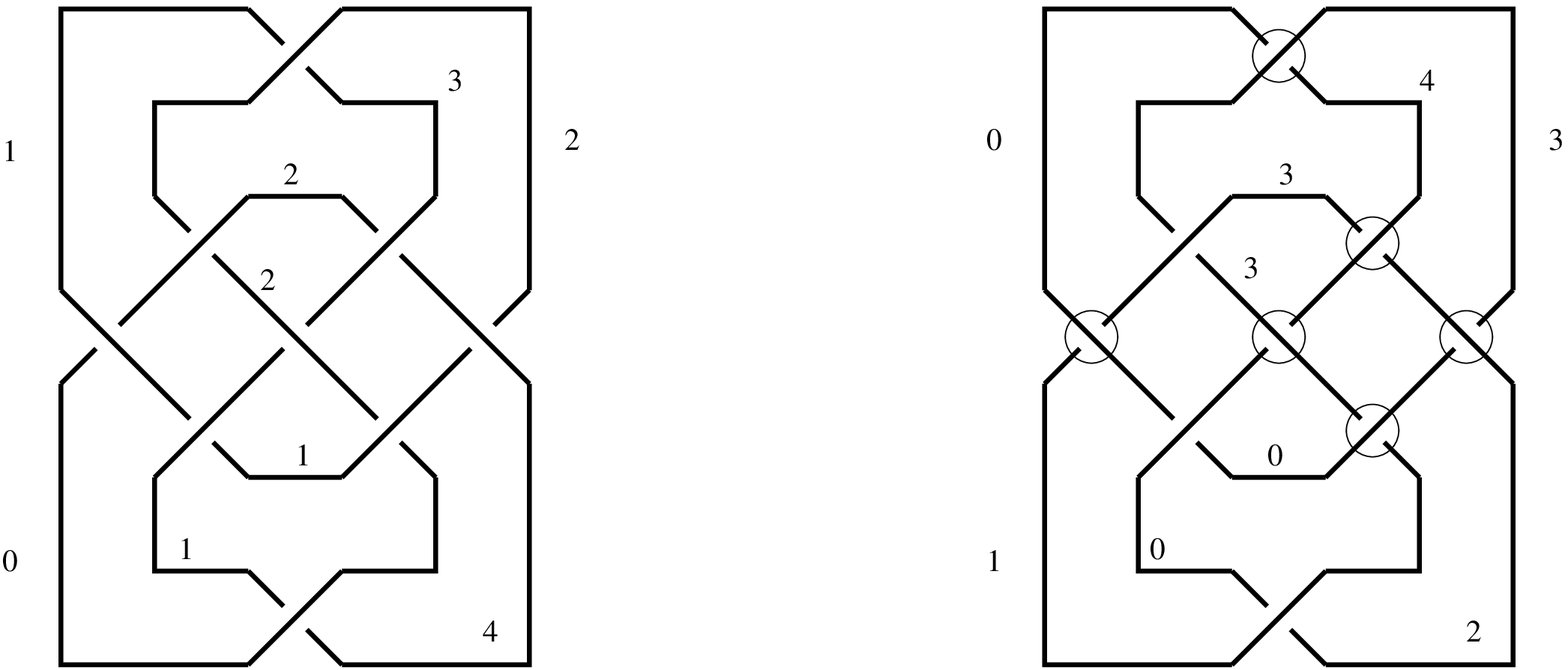}}}
	\caption{On the left, knot $9_{40}$ endowed with the $5$-coloring generated by the triplet $(0, 1, 2)$. On the right the action of permutation $(0 1) (2 3 4)$ on the colors of the coloring on the left: the result is not a $5$-coloring (at the circled crossings the coloring condition is not satisfied). }\label{fig:9_40bis}
\end{figure}

In the course of our work on colorings, we have observed that for some choices of a modulus $m>1$ and a  link admitting non-trivial $m$-colorings, the following occurs. There are distinct non-trivial $m$-colorings, $\cal C, \cal C'$ (realized on an otherwise arbitrary diagram $D$ of this link) and there is a permutation $\gamma$ of the $m$ colors such that, for each arc $a$ of $D$, the colors assigned to $a$ in the coloring $\cal C$, say ${\cal C}(a)$, and in the coloring $\cal C'$, say ${\cal C'}(a)$, satisfy:
\[
{\cal C'}(a) = \gamma \big({\cal C}(a)\big)
\]
Two such colorings will be said ``related''.  An instance where this occurs is depicted in Figure  \ref{fig:9_40tri}.
On the other hand it is not true that any permutation transforms  the colors of a coloring into the colors of another coloring (see Figure \ref{fig:9_40bis}).

\begin{figure}[!ht]
	\psfrag{a1}{\huge$c_1$}
	\psfrag{a2}{\huge$c_{2}$}
	\psfrag{a3}{\huge$c_{3}$}
    \psfrag{x}{\huge\begin{tabular}{ | c | c | c | }
\hline
$c_1$  & $c_2$  &  $c_3$ \\ \hline
0      &   0    &   1    \\ \hline
0      &   1    &   0    \\ \hline
0      &   1    &   1    \\ \hline
0      &   1    &   2    \\ \hline
0      &   1    &   3    \\ \hline
0      &   1    &   4    \\ \hline
\end{tabular}
}
	\centerline{\scalebox{.5}{\includegraphics{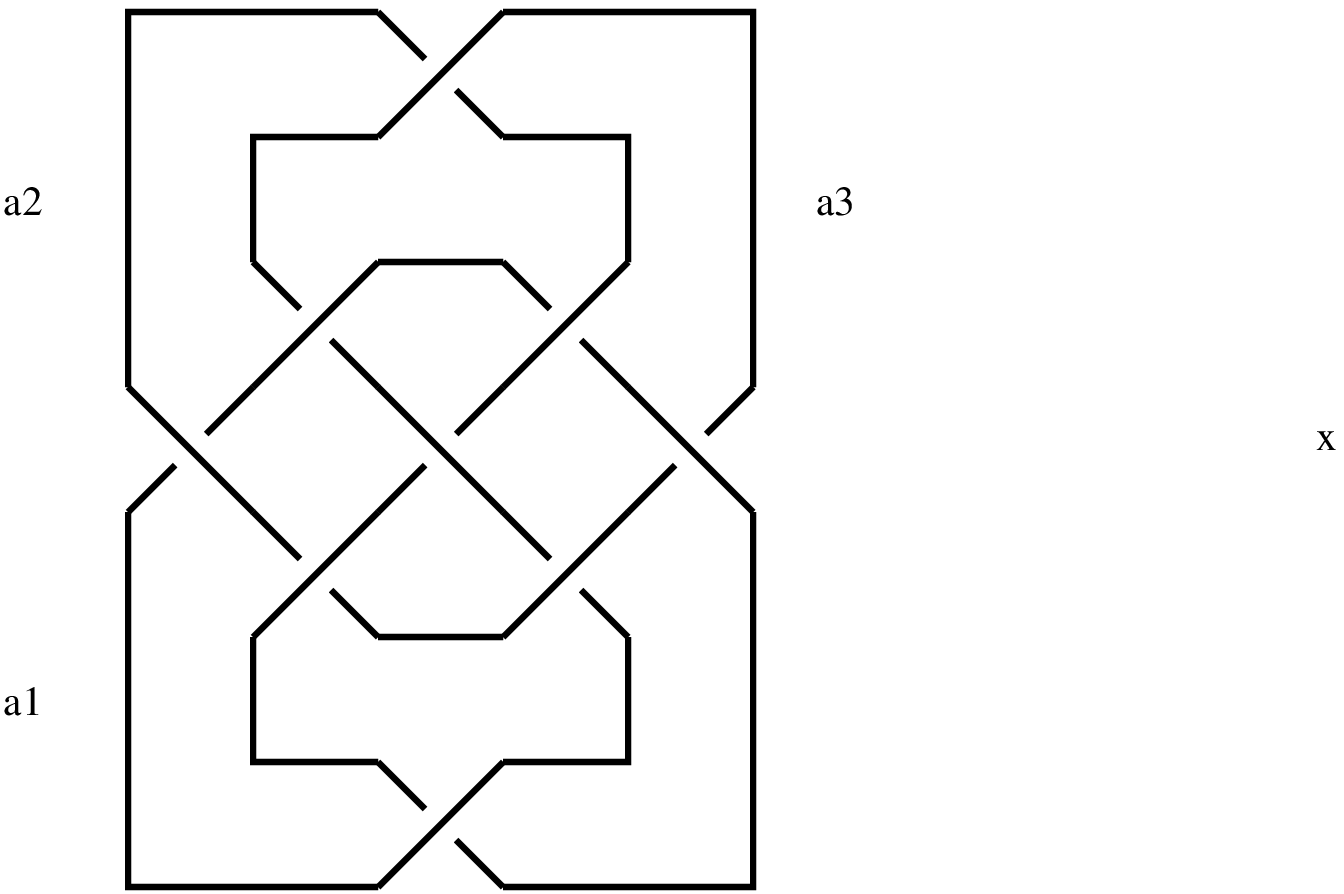}}}
	\caption{On the left, knot $9_{40}$. Assigning colors to the indicated $c_i$'s will generate a coloring of the diagram (in the sense that the other colors of this coloring are uniquely determined by $c_1, c_2, c_3$ -  we elaborate on this issue below in the text). The table on the right displays six triplets $(c_1, c_2, c_3)$ which generate essentially distinct $5$-colorings on the diagram on the left. }\label{fig:9_40}
\end{figure}

Moreover, given non-trivial $m$-colorings $\cal C$ and $\cal C'$, realized on the same diagram, it may happen that there is no permutation $\gamma$ of the $m$ colors such that for each arc $a$ of $D$
\[
{\cal C'}(a) = \gamma \big({\cal C}(a)\big)
\]
We will then say ``$\cal C'$ is essentially distinct from $\cal C$'', in the given modulus, and the colorings split into equivalence classes (to be elaborated upon below). In Figure \ref{fig:9_40} we list representatives of the distinct equivalence classes of the non-trivial $5$-colorings of $9_{40}$.

We will be primarily concerned with permutations that preserve the coloring equation at each crossing for these are the ones that actually give us a corresponding coloring of the link and we will show that the relation sketched above among $m$-colorings of a diagram is an equivalence relation (see below).

\bigbreak

We remark that the articles \cite{bot0} and \cite{bot} address the same topic as the current article. Their definition of equivalent colorings assumes one has a list of all non-trivial $m$-colorings for a given diagram and states simply that any two of these colorings are equivalent provided there is a permutation of the $m$ colors that, for each arc in the diagram, sends the color in this arc in the source coloring to the color in the same arc in the target coloring. This is equivalent to our definition. Unfortunately, for the purposes of counting equivalence classes of colorings in generic cases, the methodology in \cite{bot0} and \cite{bot} seems to resort to generating classes of colorings by letting the symmetric group on the $m$ colors act on a given $m$-coloring. As we see in Figure \ref{fig:9_40bis}, there are assignments of colors to a diagram obtained in this way that do not constitute colorings. The formulas in the articles referred to above predict in general less equivalence classes than ours due to their over-counting of the elements on each orbit.

\bigbreak

The equivalence classes of colorings constitute a topological invariant and in this article we provide combinatorial information about them. We hope this will prove to be useful for topological purposes.

\bigbreak

In Section \ref{sect:prelim} we discuss preliminaries such as the nullity and the generating arcs of a coloring (Subsection \ref{subsect:nullity}), and the definition of the equivalence classes (Subsection \ref{subsect:equivclass}). In Section \ref{sect:results} we calculate the number of equivalence classes in an infinite number of instances.

\bigbreak

\section{Preliminary Material} \label{sect:prelim}

\subsection{Nullity and Generating Arcs of a Coloring on a Diagram}\label{subsect:nullity}

Consider a link, $L$, along with a diagram $D_L$ for that link. Regarding the arcs of this diagram as algebraic variables we write the homogeneous system of linear equations consisting of the equations read off each crossing as illustrated in Figure \ref{fig:xtop}. We call the matrix of the coefficients of this homogeneous system of linear equations {\bf the coloring matrix of $D_L$}.


Any coloring matrix is made up of integers. Specifically, along each row one finds exactly two 1's and one -2, the rest being perhaps 0's. Thus, adding all the columns of a coloring matrix we obtain a column made up of 0's. It follows that the determinant of any coloring matrix is 0.

Upon performance of Reidemeister moves on a diagram, the changes on the original coloring matrix are realized by operations that constitute a subset of the following operations on integer matrices. These operations are generated by
\begin{enumerate}
\item multiplication of a row (column) by $-1$;
\item addition to one row (column) of integer linear combinations of other rows (columns);
\item insertion (deletion) of a row and column made up of 0's except for a 1 at the diagonal entry;
\item permutations of rows (columns).
\end{enumerate}
These are the operations which relate equivalent matrices over the integers (see \cite{Lickorish}, page 50). So the equivalence class of a coloring matrix is a topological invariant of the link under study. For each of these equivalence classes of matrices over the integers there is an outstanding representative which is called the Smith Normal Form (see \cite{Smith}). Although the Smith Normal Form (SNF) is a familiar object we elaborate here slightly about it in order to bring out some connections with colorings of knots which we do not find in the literature.

An integer matrix in Smith Normal Form is a matrix such that its entries are all zero except perhaps along the diagonal. Along the diagonal the entries are non-negative (without loss of generality) and the $i$-th entry divides the $(i+1)$-th entry, up to a certain index $l$, and after that, the entries are all $0$'s:
\[
d_1, d_2, \dots , d_l, 0, 0, \dots , 0 \qquad \text{ with } \qquad d_i | d_{i+1} \, \quad 1\leq i \leq l-1
\]

The $d_i$'s are called the {\bf invariant factors} of the equivalence class; their name reflects the fact that the multi-set formed by them is an invariant of the equivalence class. This multi-set is then a topological invariant if it originates from a coloring matrix. Moreover, the Smith Normal Form of a coloring matrix is sure to have a $0$ at the last entry of the diagonal since we proved above that the determinant of a coloring matrix is $0$. The product of the remaining entries of the diagonal of the Smith Normal Form of a coloring matrix is {\bf the determinant of the link} under study. This is also a topological invariant.  (In passing, it is known that for knots i.e., $1$-component links, the determinant of the knot is an odd integer, see \cite{Rolfsen}.)

We denote the Smith Normal Form of a matrix $M$ by $S(M)$. Being an element of the equivalence class of $M$, $S(M)$ is obtained by a finite number of the operations listed above. We may then collect all the information concerning the row operations into an invertible matrix called $R$ and likewise for the column operations into an invertible matrix called $C$ to state (\cite{Smith})
\begin{align}\label{eqn:eqn}
S(M) = R  M C
\end{align}
with the juxtaposition of pairs of consecutive symbols on the right-hand side of the equation denoting matrix multiplication.

Let us now fix an otherwise arbitrary link along with one of its diagrams.
Let us then relate the Smith Normal Form (and its invariant factors) of the coloring matrix of this diagram to the corresponding system of linear homogeneous equations and its solutions. There are always solutions of this system of equations namely by assigning the same integer to each arc. This corresponds to the fact that the determinant of the coloring matrix is $0$. One of the algebraic variables may take on any value and if there is no other zero entry along the diagonal of the Smith Normal Form, then the remaining variables are uniquely determined once the former variable has been assigned a value. Going back to the original system of equations we obtain the so-called trivial solutions i.e., those solutions that assign the same value to each and every arc of the diagram.

The invariant factors associated to our coloring matrix via its Smith Normal Form allow us to do something else. Suppose we choose a factor $m$ of one of these invariant factors and decide to work over the integers modulo $m$. Then our Smith Normal Form in this new setting acquired at least one more $0$ along the diagonal. Then, there is at least one more variable which can take on any value, modulo $m$. Going back to the original system of equations, there are at least two arcs which can take on any value modulo $m$. Hence, we now have polychromatic colorings i.e., solutions where at least two distinct arcs take on two distinct colors that is, values modulo $m$. Had we chosen an $m$ which does not possess common factors with the invariant factors, then modulo $m$ there would have been only trivial colorings.

\begin{prop}\label{prop:genarcsp}Let $p$ be an odd prime. Let $D$ be a link diagram.
The number of $0$'s $($modulo $p)$ along the diagonal of the Smith Normal Form of the coloring matrix, $M$, of $D$, equals the least number of  arcs of $D$ that can independently receive colors modulo $p$, and generate each $p$-coloring of $D$.
\end{prop}Proof.
If the Smith Normal Form exhibits $n$ $0$'s modulo  $p$, this means that the space of solutions has dimension $n$; working modulo a prime implies we are doing Linear Algebra over a field so it makes sense to talk about dimensions of spaces and bases. Then matrix $C$ in (\ref{eqn:eqn}) above operates a change of basis taking us back to algebraic variables equivalent to the arcs of the original diagram. Then $n$ of these arcs have to generate all the colorings (i.e., all the solutions of the indicated system of equations modulo $p$) in terms of a basis of coloring vectors.
$\hfill \blacksquare$

\begin{def.}\label{def:pnul}
The number $n$ in the proof of Proposition \ref{prop:genarcsp}  is called the $p$-nullity $($or the rank mod $p )$ of the coloring matrix of $D$. Any set of  $n$ arcs that can independently receive colors modulo $p$ and so generate each $p$-coloring of the diagram under study is said a set of generating arcs $($of this diagram, with respect to this modulus$)$.
\end{def.}

\begin{cor}
We keep the notation of Proposition \ref{prop:genarcsp}. If the $p$-nullity of a link is $n$ then there are $p^n$ $p$-colorings of the link, and $p^n-p$ non-trivial $p$-colorings of this link.
\end{cor}Proof. There  are $p$ integers mod $p$ so there are always $p$ trivial $p$-colorings and $p^n$ $p$-colorings.
$\hfill \blacksquare$

\begin{cor}
Let  $m$ be a composite positive integer. Let $D$ be a link diagram. Each zero $($modulo $m)$ along the diagonal of the Smith Normal Form of the coloring matrix, $M$, of $D$ contributes with a factor $m$ for the number of solutions. Each zero divisor, $z$, of $m$ along the diagonal contributes with a factor $\gcd (z, m)$ to the number of solutions. With $n_Z$ for the number of $0$'s $($modulo $m)$ along the diagonal in $S(M)$, and $IZ(M)$ for the set the of invariant factors of $M$ which are zero divisors of $m$, the formula for the number of $m$-colorings of $D$ is:
\[
m^{n_Z}\prod_{z\in IZ(M)}\gcd (z, m)
\]
\end{cor}Proof. The contribution of the $n_Z$ zero's (modulo $m$) along the diagonal of the Smith Normal Form to the number of solutions is clear. For the contribution of the zero divisors along the diagonal to the number of solutions see  \cite{kl}, page 40. This concludes the proof.
$\hfill \blacksquare$

\bigbreak

\subsection{Equivalence Classes of Colorings}\label{subsect:equivclass}

\noindent

In this section we introduce equivalence classes of colorings as orbits of actions of certain groups of permutations on the set of colorings of a diagram. In order for this notion to be topological we require a special kind of permutation which we call a {\bf Coloring Automorphism}. These are permutations which comply with the coloring operation,
\[
a \ast b := 2b-a
\]
in a pre-assigned modulus $m$. This operation generalizes to the quandle operation, generalizing also the notion of coloring (\cite{dJoyce, m}). In the particular instance $a \ast b = 2b-a$ we are dealing with the so-called dihedral quandles, one per integer modulus $m$.

\begin{def.}[Coloring Automorphism of $\mathbf{Z}_m$]\label{def:auto}
Given an integer $m\geq 3$, we define a coloring automorphism of $\mathbf{Z}_m$ to be a permutation, $f$, of $\mathbf{Z}_m$ such that
\[
f(a\ast b) = f(a) \ast f(b)
\]
for all $a, b \in \mathbf{Z}_m$, with $x\ast y = 2y-x$ $($mod $m )$, for every $x, y\in \mathbf{Z}_m$.
\end{def.}

In \cite{elhamdadi} we find the following facts. For a given integer $m\geq 3$, each coloring automorphism of $\mathbf{Z}_m$ is given by:
\[
f_{\lambda, \mu}(x) = \lambda x + \mu
\]
with $\mu \in \mathbf{Z}_m$ and $\lambda \in \mathbf{Z}_m^{\ast}$, the set of units of $\mathbf{Z}_m$. The set of all these coloring automorphisms of $\mathbf{Z}_m$ equipped with composition of functions, constitutes a group isomorphic to the affine group over $\mathbf{Z}_m$ i.e., isomorphic to the semi-direct product $\mathbf{Z}_m \rtimes\mathbf{Z}_m^{\ast}$. We denote it $\mathbf{Aut}_m$.

For any integer $m\geq 3$ the inner coloring automorphism group of $\mathbf{Z}_m$  is generated by the automorphisms of the form $f_b(x)=x\ast b$. It is easy to see that this group consists of the elements of the form,
\[
f_{\pm, \mu}(x)=\pm x + \mu
\]
If $m$ is even this subgroup is isomorphic to the dihedral group of order $m$ and $\mu$ can take on only ``even'' values from $\mathbf{Z}_m$. If $m$ is odd, this subgroup is isomorphic to the dihedral group of order $2m$ and $\mu$ can take on any value from $\mathbf{Z}_m$. We denote it $\mathbf{Inn}_m$. This information about coloring automorphisms of $\mathbf{Z}_m$ is contained in \cite{elhamdadi}.

\bigbreak

In the sequel, we will write ``automorphism'' (respectively, ``inner automorphism'') instead of  the longer ``coloring automorphism of $\mathbf{Z}_m$'' (respectively, ``inner coloring automorphism of $\mathbf{Z}_m$'') since these are the only automorphisms of $\mathbf{Z}_m$ we consider in this article i.e., the permutations of elements of $\mathbf{Z}_m$ that comply with the coloring operation.

Specifically, we will use the expression {\bf automorphism} to designate a permutation of the form

\[
f_{\lambda, \mu}(x) = \lambda x + \mu
\]
with $\mu \in \mathbf{Z}_m$ and $\lambda \in \mathbf{Z}_m^{\ast}$, and {\bf inner automorphism} to designate a permutation of the form
\[
f_{\pm, \mu}(x)=\pm x + \mu
\]
with $\mu$ taking on only ``even'' values from $\mathbf{Z}_m$ if $m$ is even; with $\mu$ taking on any value from $\mathbf{Z}_m$ if $m$ is odd.

\bigbreak
We remark that it is well known that for a quandle $(Q,*)$ and
a diagram $D$, the set of diagram
colorings by elements of $Q$, $Col_Q(D)$ is a Q-quandle set, where the
action of $Q$ on $Col_Q(D)$ is given
by $C*q$ for a coloring $C\in Col_Q(D)$ and $q\in Q$ (Kamada was
the first proponent of this language). Our considerations for dihedral quandles are related to this.

\bigbreak
\begin{def.} Let $m > 1 $ be an integer. Let $L$ be a link admitting non-trivial $m$-colorings. Let $D$  be a diagram of $L$. We let $\mathbf{m{\cal C}D}$ stand for the set of {\bf non-trivial} $m$-colorings of $D$.
\end{def.}
\bigbreak
\begin{prop}\label{prop:action} Let $m > 1 $ be an integer. Let $L$ be a link admitting non-trivial $m$-colorings and let $D$ be a diagram of $L$. Let $G$ be a subgroup of $Aut_m$. Then $G$ acts on $m{\cal C}D$ by permutations.

Specifically, given $g\in G$ and ${\cal C}$, an $m$-coloring of $D$ with colors $c_i$, then  $g{\cal C}$ is the $m$-coloring of $D$ obtained by replacing each color $c_i$ by $g(c_i)$.

Moreover, this action is faithful and if $m$ is prime this action is also free.
\end{prop}{\bf Proof}. We keep the notation of the statement.  We regard ${\cal C}\in m{\cal C}D$ as the map which assigns colors to the arcs of $D$ in such a way that, ${\cal C}(a_{i+1})=2{\cal C}(a_{j_i})-{\cal C}(a_i)$, where $j_i$ designates the index of the over-arc of the crossing where under-arcs with indices $i$ and $i+1$ meet, see Figure \ref{fig:xtop} (where now each $c_k$ should be read ${\cal C}(a_k)$).

So, given $g\in G$ and ${\cal C}\in m{\cal C}D$, then $g{\cal C}$ is such that
\[
g\big({\cal C}(a_{i+1})\big)=g\big(2{\cal C}(a_{j_i})-{\cal C}(a_i)\big)=2g\big({\cal C}(a_{j_i})\big)-g\big({\cal C}(a_i)\big)
\]
so $g{\cal C}$ is again an $m$-coloring of $D$.

Clearly, the identity element $1_G\in G$ is such that $1_G {\cal C} = {\cal C}$. Furthermore, for any two $g_1, g_2\in G$, the composition of functions guarantees that $(g_1g_2){\cal C} = g_1\big(g_2{\cal C}\big)$.

We now prove that this action is faithful i.e., we prove that given a non-identity $g\in G$ there exists a coloring ${\cal C}\in m{\cal C}D$ such that $g{\cal C}\in m{\cal C}D \neq {\cal C}\in m{\cal C}D$. We recall that the elements of $G$ are, in particular, permutations of the elements of $\mathbf{Z}_m$. So given a non-identity element of $G$ which moves $i\in \mathbf{Z}_m$, then the coloring obtained by assigning $i$ to one of the generating arcs of the diagram is transformed via $g$ into a coloring where now this generating arc is assigned $g(i)\neq i$.

We now prove that this action is free i.e., that if given $g, h\in G$ there exists a coloring ${\cal C}\in m{\cal C}D$ such that $g\big({\cal C}\big) =h\big({\cal C}\big)$ then $g=h$. We recall that, for some $\lambda , \lambda' \in \mathbf{Z}_m^{\ast}$ and $\mu , \mu' \in \mathbf{Z}_m$, $g (x)=\lambda x + \mu$, $h(x)=\lambda' x+\mu'$ for any $x\in \mathbf{Z}_m$. Since $g\big({\cal C}\big) =h\big({\cal C}\big)$ then there exists two distinct colors in $\mathbf{Z}_m$, say $a\neq b$ such that $g(a)=h(a)$ and $g(b)=h(b)$. More precisely,
\begin{equation*}
\begin{cases}
0=(\lambda -\lambda')a + (\mu -\mu') \\
0=(\lambda -\lambda')b + (\mu -\mu')
\end{cases}
\quad \Longleftrightarrow \quad
\begin{cases}
\lambda = \lambda' \\
\mu = \mu'
\end{cases}
\end{equation*}
since $a \neq b$ and $m$ is prime.  Thus $g=h$.

This concludes the proof.

$\hfill \blacksquare$

\begin{def.}[$G$-Equivalence Classes of $m$-Colorings of $D$]\label{def:GequivD}
Let $G$ be a subgroup of $\mathbf{Aut}_m$.

The $G$-Equivalence Classes of $m$-Colorings of $D$ are, by definition, the $G$-orbits over $m{\cal C}D$.
\end{def.}

We will next prove that this notion provides topological invariants (in particular, the number of equivalence classes of $m$-colorings of a link).

It is well known that, given any two diagrams of the same link, there is a  bijection between the two sets of $m$-colorings of these diagrams (\cite{pLopes}, \cite{Przytycki}, Lemma $2.2$). Moreover, this bijection takes trivial $m$-colorings to trivial $m$-colorings and non-trivial $m$-colorings to non-trivial $m$-colorings. This bijection is realized by the ``Colored Reidemeister Moves''. The Colored Reidemeister Moves apply to a diagram endowed with an $m$-coloring; a Reidemeister move is applied to the diagram and a local adjustment of the coloring is performed. These adjustements are unique and reversible thereby proving the bijection between the two sets of $m$-colorings of any two diagrams of the same link.

\bigbreak

\begin{figure}[!ht]
	\psfrag{D}{\huge$D$}
	\psfrag{D'}{\huge$D'$}
	\psfrag{calCr=x}{\huge${\cal C}(r)=x$}
	\psfrag{fcalCr=x}{\huge$g\big( {\cal C}(r)\big)=g(x)$}
	\psfrag{calC'r1r2=x}{\huge${\cal C'}(r'_1)=x={\cal C'}(r'_2)$}
	\psfrag{fcalC'r1r2=x}{\huge$g\big( {\cal C'}(r'_1)\big)=g(x)=g\big( {\cal C'}(r'_2)\big)$}
	\psfrag{r}{\huge$r$}
	\psfrag{r1}{\huge$r'_1$}
	\psfrag{r2}{\huge$r'_2$}
	\psfrag{tilde}{\huge$\sim$}
	\psfrag{9}{\huge$\mathbf{9}$}
	\centerline{\scalebox{.5}{\includegraphics{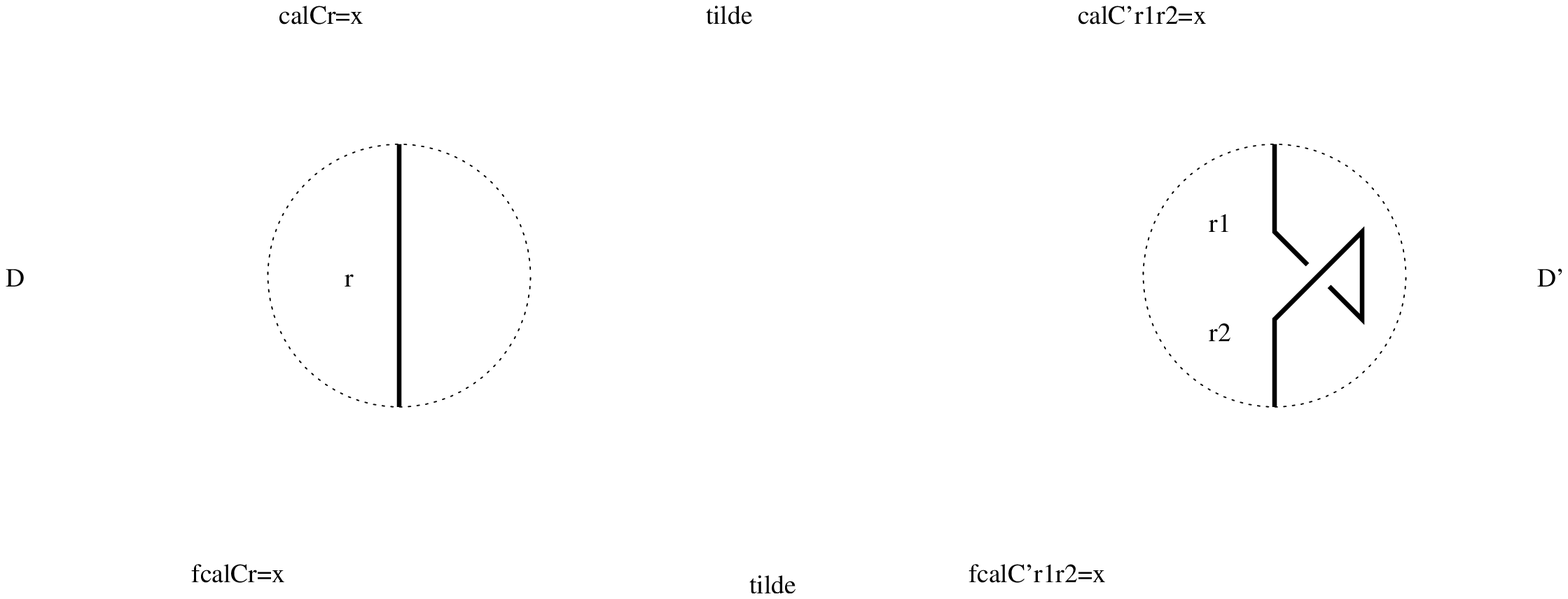}}}
	\caption{Colored Reidemeister move of type I and $G$-equivalence relation of colorings on the same diagram.}\label{fig:r1}
\end{figure}

\begin{figure}[!ht]
	\psfrag{D}{\huge$D$}
	\psfrag{D'}{\huge$D'$}
	\psfrag{calCr=x}{\huge${\cal C}(r)=x$}
	\psfrag{calCs=y}{\huge${\cal C}(s)=y$}
	\psfrag{fcalCr=x}{\huge$g\big( {\cal C}(r)\big)=g(x)$}
	\psfrag{fcalCs=y}{\huge$g\big( {\cal C}(s)\big)=g(y)$}
	\psfrag{calC'r1=x}{\huge${\cal C'}(r'_1)=x \qquad {\cal C'}(s')=y$}
	\psfrag{calC'r2=2y-x}{\huge${\cal C'}(r'_2)=2{\cal C'}(s') - {\cal C'}(r'_1) = 2y-x$}
	\psfrag{calC'r3=x}{\huge${\cal C'}(r'_3)= \cdots = 2y- (2y-x) = x$}
	\psfrag{fcalC'r1=x}{\huge$g\big( {\cal C'}(r'_1)\big)=g(x) \qquad g\big( {\cal C'}(s')\big)=g(y)$}
	\psfrag{fcalC'r2=2y-x}{\huge${g\big( \cal C'}(r'_2)\big)=g(2{\cal C'}(s') - {\cal C'}(r'_1)) = 2g(y)-g(x)$}
	\psfrag{fcalC'r3=x}{\huge${g\big( \cal C'}(r'_3)\big)= \cdots = g(2y- (2y-x)) = g(x)$}
	\psfrag{r}{\huge$r$}
	\psfrag{s}{\huge$s$}
	\psfrag{s'}{\huge$s'$}
	\psfrag{r1}{\huge$r'_1$}
	\psfrag{r2}{\huge$r'_2$}
	\psfrag{r3}{\huge$r'_3$}
	\psfrag{tilde}{\huge$\sim$}
	\psfrag{9}{\huge$\mathbf{9}$}
	\centerline{\scalebox{.5}{\includegraphics{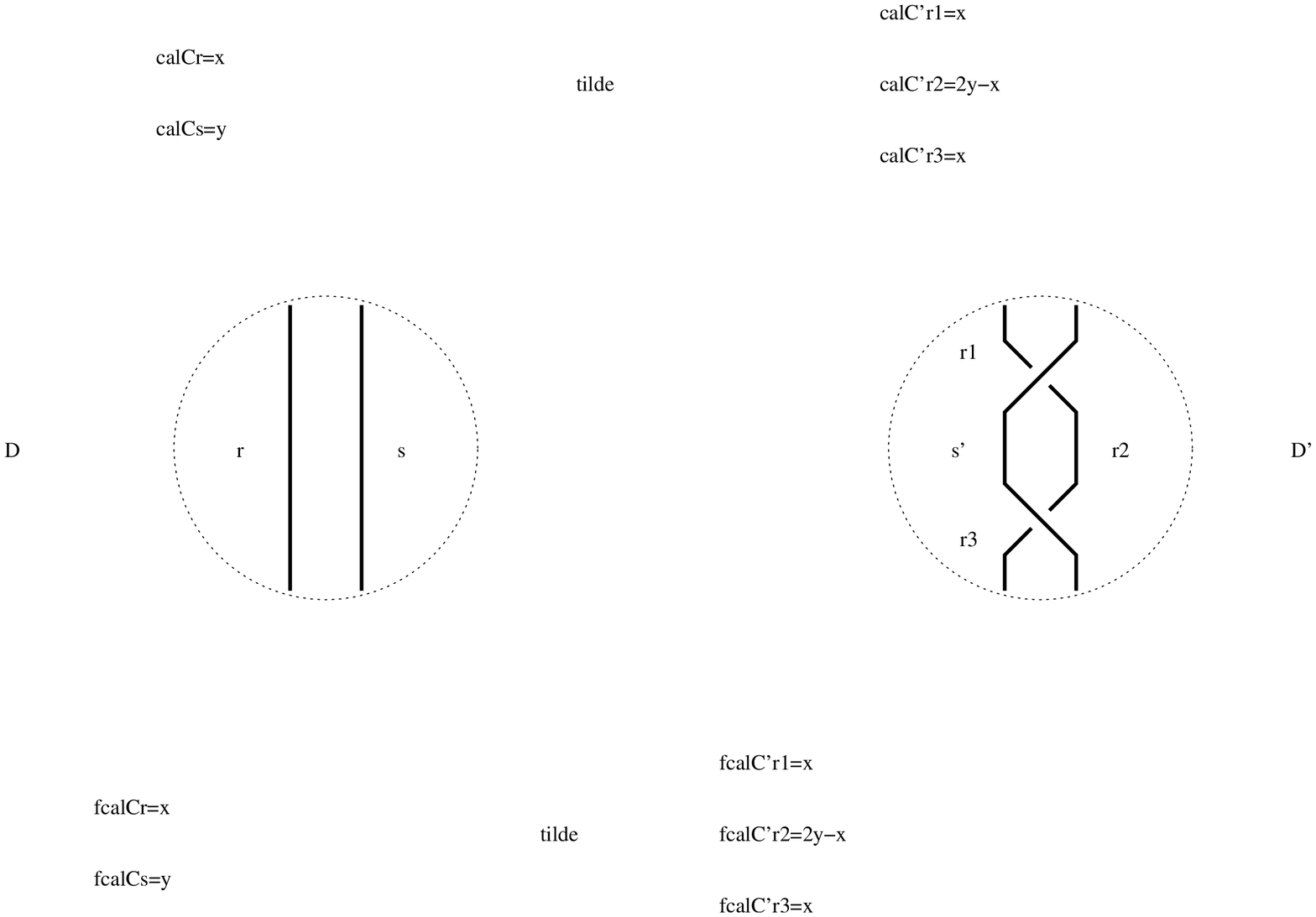}}}
	\caption{Colored Reidemeister move of type II and $G$-equivalence relation of colorings on the same diagram.}\label{fig:r2}
\end{figure}

\begin{figure}[!ht]
	\psfrag{D}{\huge$D$}
	\psfrag{D'}{\huge$D'$}
	\psfrag{calCr1}{\huge${\cal C}(r_1)=x\qquad {\cal C}(s_1)=y\qquad {\cal C}(t)=z$}
	\psfrag{calCr2}{\huge${\cal C}(r_2)=\cdots =2y-x$}
	\psfrag{calCs2}{\huge${\cal C}(s_2)=\cdots =2z-y$}
	\psfrag{calCr3}{\huge${\cal C}(r_3)=\cdots =2z-(2y-x)=2z-2y+x$}
	\psfrag{calC'r1}{\huge${\cal C'}(r'_1)=x\qquad {\cal C'}(s'_1)=y\qquad {\cal C'}(t')=z$}
	\psfrag{calC'r2}{\huge${\cal C'}(r'_2)=\cdots =2z-x$}
	\psfrag{calC's2}{\huge${\cal C'}(s'_2)=\cdots =2z-y$}
	\psfrag{calC'r3}{\huge${\cal C'}(r'_3)=\cdots =2z-2y+x$}
	\psfrag{fcalCr1}{\huge$g\big( {\cal C}(r_1)\big)=g(x)\qquad g\big( {\cal C}(s_1)\big)=g(y)\qquad g\big( {\cal C}(t)\big)=g(z)$}
	\psfrag{fcalCr2}{\huge$g\big( {\cal C}(r_2)\big)=\cdots =2g(y)-g(x)$}
	\psfrag{fcalCs2}{\huge$g\big( {\cal C}(s_2)\big)=\cdots =2g(z)-g(y)$}
	\psfrag{fcalCr3}{\huge$g\big( {\cal C}(r_3)\big)=\cdots =2g(z)-2g(y)+g(x)$}
	\psfrag{fcalC'r1}{\huge$g\big( {\cal C'}(r'_1)\big)=g(x)\quad g\big( {\cal C}(s_1)\big)=g(y)\quad g\big( {\cal C}(t)\big)=g(z)$}
	\psfrag{fcalC'r2}{\huge$g\big( {\cal C'}(r'_2)\big)=\cdots =2g(y)-g(x)$}
	\psfrag{fcalC's2}{\huge$g\big( {\cal C'}(s'_2)\big)=\cdots =2g(z)-g(y)$}
	\psfrag{fcalC'r3}{\huge$g\big( {\cal C'}(r'_3)\big)=\cdots =2g(z)-2g(y)+g(x)$}
	\psfrag{t}{\huge$t$}
	\psfrag{t'}{\huge$t'$}
	\psfrag{r1}{\huge$r_1$}
	\psfrag{r2}{\huge$r_2$}
	\psfrag{r3}{\huge$r_3$}
	\psfrag{r'1}{\huge$r'_1$}
	\psfrag{r'2}{\huge$r'_2$}
	\psfrag{r'3}{\huge$r'_3$}
	\psfrag{s1}{\huge$s_1$}
	\psfrag{s'1}{\huge$s'_1$}
	\psfrag{s2}{\huge$s_2$}
	\psfrag{s'2}{\huge$s'_2$}
	\psfrag{tilde}{\huge$\sim$}
	\psfrag{9}{\huge$\mathbf{9}$}
	\centerline{\scalebox{.35}{\includegraphics{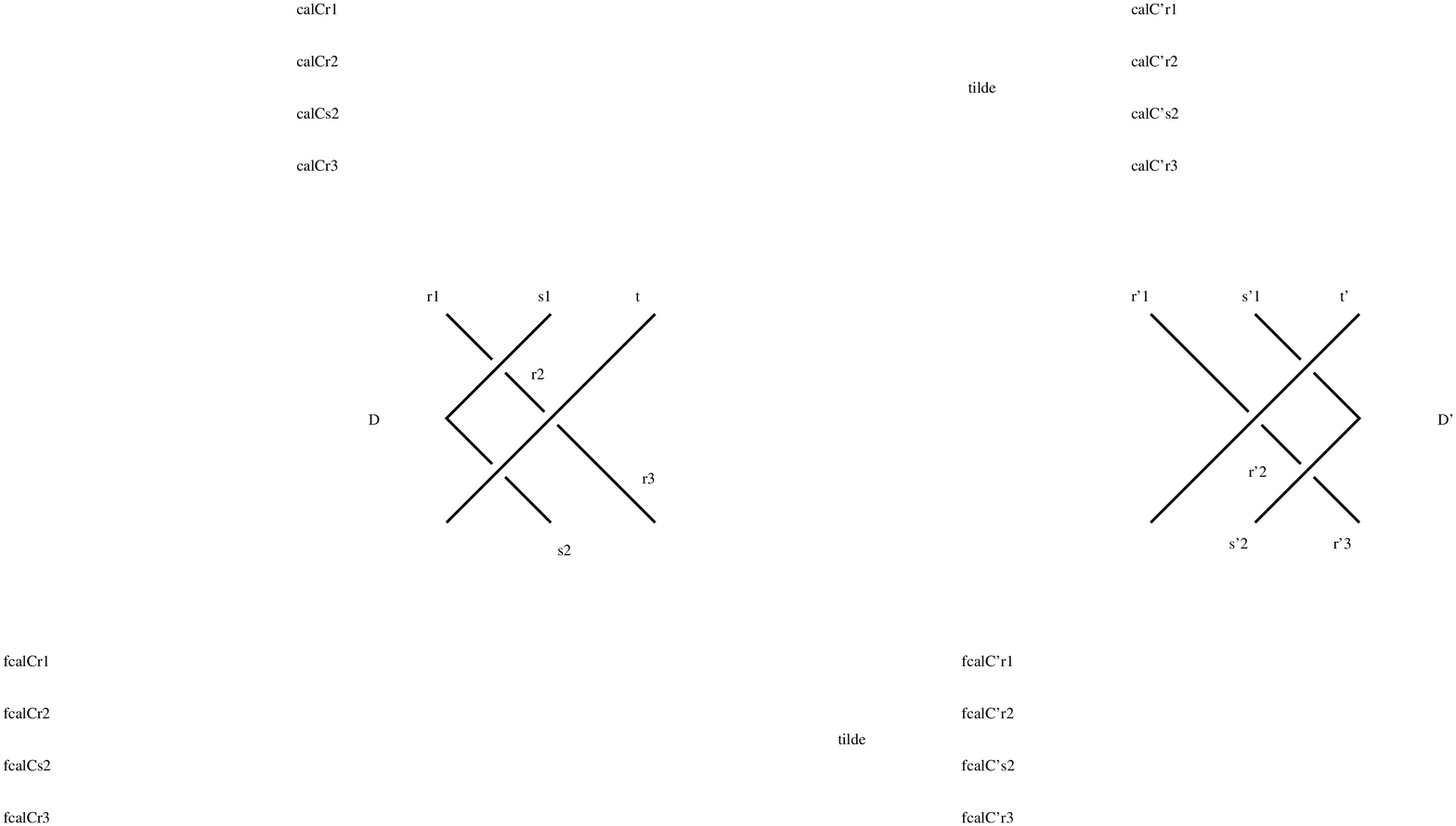}}}
	\caption{Colored Reidemeister move of type III and $G$-equivalence relation of colorings on the same diagram.}\label{fig:r3}
\end{figure}

\begin{prop}\label{prop:equivtop}
Let $m$ be an integer greater than $1$. Let $L$ be a link admitting non-trivial $m$-colorings and let  $D$ and $D'$ be two diagrams of $L$. Let $G$ be a subgroup of $\mathbf{Aut}_m$.

There is a bijection from $m{\cal C}D$ to $m{\cal C}D'$, which preserves the $G$-equivalence classes.
\end{prop}{\bf Proof}: From \cite{pLopes} we know that the Colored Reidemeister Moves realize a bijection from the set of $m$-colorings of $D$ to the set of $m$ colorings of $D'$, taking non-trivial colorings to non-trivial colorings. We now prove that the Colored Reidemeister Moves take distinct elements of $m{\cal C}D$ along a $G$-equivalence class, to distinct elements of $m{\cal C}D'$ along a $G$-equivalence class. Specifically, for $g\in G$ and ${\cal C}\in m{\cal C}D$, we prove that the ``Colored Reidemeister moves'' take ${\cal C}\in m{\cal C}D$ to ${\cal C'}\in m{\cal C}D'$ and $g{\cal C}\in m{\cal C}D$ to $g{\cal C'}\in m{\cal C}D'$. The proofs of these statements for the individual ``Colored Reidemeister Moves'' of type I, II, and III are displayed in Figures \ref{fig:r1}, \ref{fig:r2}, and \ref{fig:r3}. The $\sim$ associates horizontally colorings on distinct diagrams related by a Colored Reidemeister move. Vertically we display colorings ${\cal C}$ and $g{\cal C}$ (${\cal C'}$ and $g{\cal C'}$, respect.) for diagram $D$ ($D'$, respect.).

In Figures \ref{fig:r1} and \ref{fig:r2} circles with dotted lines were drawn to bring out the local nature of the transformation. This was not done in Figure \ref{fig:r3} in order not to overburden the Figure.
$\hfill \blacksquare$

\bigbreak

{\bf Remark.} Proposition 2.3 can also be seen by regarding mod-m colorings as
elements of
$H_1(M^{(2)}_L),Z/m)\cong Hom(H_1((M^{(2)}_L), Z/m)$.
\bigbreak

\bigbreak

\begin{thm}\label{cor:numbers}
Let $L$ be a link and $D$ one of its diagrams. Let $m>1$ be an integer.

The number of $G$-equivalence classes of $m$-colorings of $D$ is a topological invariant. The multi-set whose elements are the number of $m$-colorings per $G$-equivalence class of $m$-colorings of $D$ is a topological invariant.
\end{thm}{\bf Proof}. This is a straight-forward consequence of Proposition \ref{prop:equivtop}.
$\hfill \blacksquare$

\bigbreak

\begin{figure}[!ht]
	\psfrag{5}{\huge$\mathbf{5}$}
	\psfrag{a=0}{\huge$a=0$}
	\psfrag{b=1}{\huge$b=1$}
	\psfrag{1}{\huge$1$}
	\psfrag{b1=2}{\huge$b=2$}
	\psfrag{2}{\huge$2$}
	\psfrag{3}{\huge$3$}
	\psfrag{4}{\huge$4$}
	\centerline{\scalebox{.5}{\includegraphics{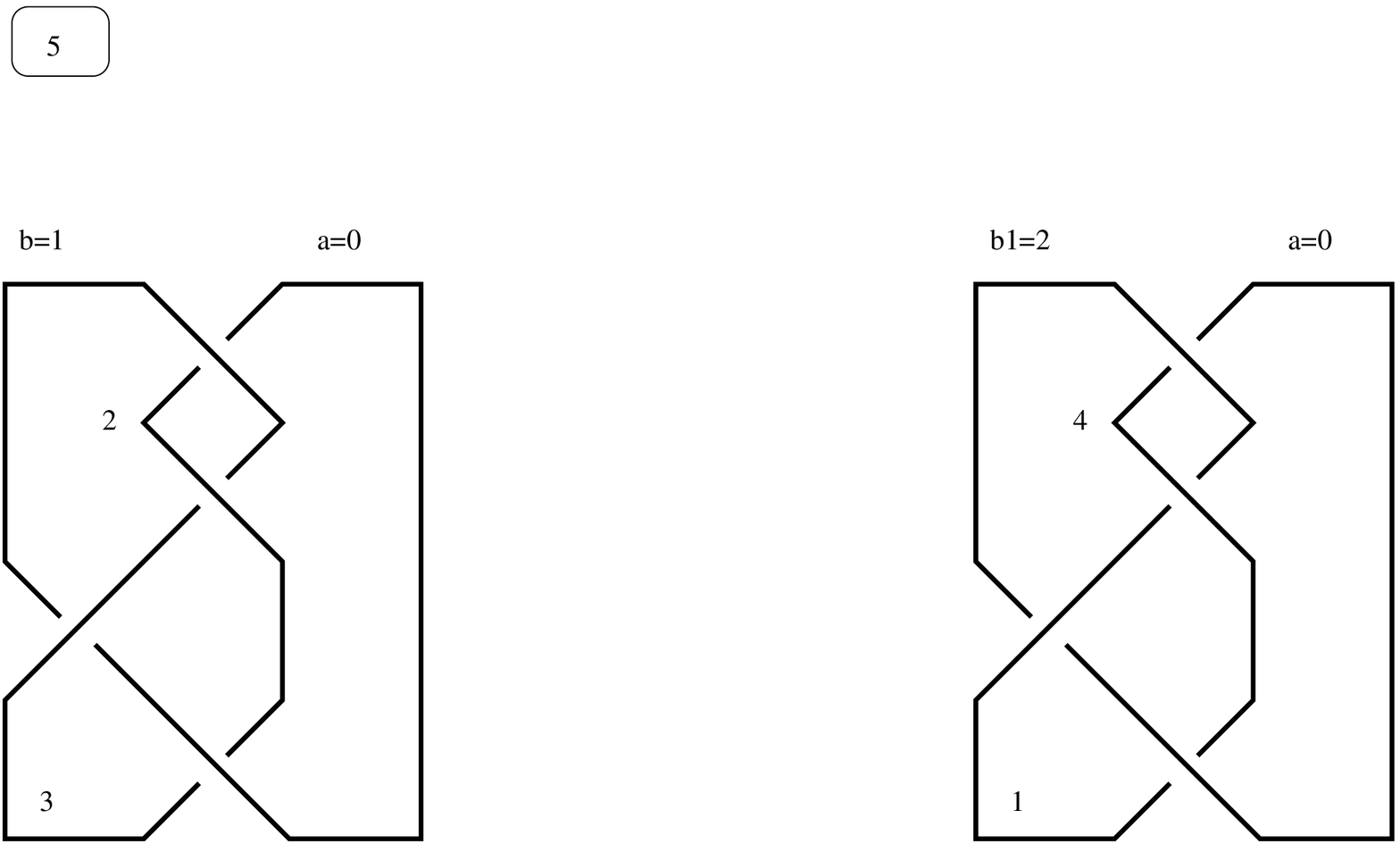}}}
	\caption{Two $5$-colorings of the Figure-$8$ knot (which has determinant $5$). These colorings are representatives of the two distinct $5$-coloring inner equivalence classes. On the other hand there is only one $5$-coloring equivalence class. These facts will be clear from the results in Section \ref{sect:results}.}\label{fig:fig8}
\end{figure}

We remark that in the sequel $G$, the subgroup of $\mathbf{Aut}_m$,  will be either $\mathbf{Aut}_m$ itself or $\mathbf{Inn}_m$. Figure \ref{fig:fig8} illustrates the fact that in general there are more inner equivalence classes than equivalence classes (for the same link and for the same modulus).

\bigbreak

\section{Formulas for Numbers of Equivalence Classes} \label{sect:results}
\noindent

In this Section we apply the theory developed above to specific situations. {\bf We use $p$-nullity as in Definition \ref{def:pnul}.}

\subsection{Equivalence Classes} \label{subsect:equivclasses}
\noindent

\begin{prop}\label{prop:numberequivclasscol}Let $p$ be an odd prime and $n$ an integer greater than $1$.
A link $L$ with $p$-nullity $n$ has
\[
\frac{p^{n-1}-1}{p-1}
\]
equivalence classes of $p$-colorings.
\end{prop}{\bf Proof}. As discussed right after Definition \ref{def:auto}, an automorphism of $\mathbf{Z}_p$
\[
f_{\lambda, \mu} (x) = \lambda x + \mu
\]
depends on two parameters $\lambda \in \mathbf{Z}_p^\ast$ and $\mu\in \mathbf{Z}_p$. Since $|\mathbf{Z}_p|=p$ and $|\mathbf{Z}_p^\ast|=p-1$, there are then exactly $p(p-1)$ automorphisms for $\mathbf{Z}_p$.

Now suppose $\cal C$ is in $p{\cal C}D$, where $D$ is a diagram of $L$. Since the action of $\mathbf{Aut}_p$ is free (\ref{prop:action}) each orbit of the action has exactly $p(p-1)$ elements. Since there are $p^n-p$ elements in $p{\cal C}D$, there are then
\[
\frac{p^n-p}{p(p-1)} = \frac{p^{n-1}-1}{p-1}
\]
orbits of this action which is the number of equivalence classes of $p$-colorings for link $L$.
$\hfill \blacksquare$

\bigbreak

\begin{cor}\label{cor:mincol}
We keep the notation of Proposition \ref{prop:numberequivclasscol}.
\begin{enumerate}
\item If a diagram $D$ of link $L$ admits a non-trivial $p$-coloring with the least number of colors (over all diagrams, over all non-trivial $p$-colorings), then there are at least $p(p-1)$ such $p$-colorings of $D$.

\item If the nullity of $L$ mod $p$ is $2$ and a diagram $D$ of $L$ admits a non-trivial $p$-coloring with $k$ colors, then any other non-trivial $p$-coloring of $D$ uses $k$ colors. In particular, if $D$ is a diagram of $L$ where a non-trivial $p$-coloring is realized with the least number of colors, then any other non-trivial $p$-coloring of this diagram uses also the least number of colors.
\end{enumerate}
\end{cor}Proof.
\begin{enumerate}
\item Since the automorphisms are permutations of the $p$ colors they preserve the number of distinct colors. So if a non-trivial $p$-coloring of a diagram uses $k$ colors, then along its equivalence class the non-trivial $p$-colorings use $k$ colors each and there are $p(p-1)$ non-trivial $p$-colorings per equivalence class. If a diagram $D$ of link $L$ admits a non-trivial $p$-coloring with the least number of colors then along its equivalence class the non-trivial $p$-colorings use the same number of colors each.

\item If the nullity of $L$ mod $p$ is $2$ then there is only $1$ equivalence class (mod $p$). Then the arguing of $1.$ is valid for the $p(p-1)$ non-trivial $p$-colorings in this orbit.
\end{enumerate}
$\hfill \blacksquare$

\bigbreak

\begin{cor}\label{cor:D}
Let $L$ be a link with the following property.

The Smith Normal Form of a$($ny$)$ coloring matrix of $L$ has only one $0$ and only one $0 \neq d \neq 1$ along the diagonal.

Then for any prime $p$ such that $p|d$, there is only one equivalence class of $p$-colorings. In particular, rational  links satisfy this property.
\end{cor}Proof. Working mod $p$ the Smith Normal Form of the coloring matrix will exhibit exactly two zeros. Hence the $p$-nullity is $2$ and the result follows from Proposition \ref{prop:numberequivclasscol}.
$\hfill \blacksquare$

\bigbreak

\begin{cor}\label{cor:otherD}
The following links have only one class of $p$-colorings for each prime $p$ for which they admit non-trivial $p$-colorings.
\begin{enumerate}
\item Links whose determinant is prime.
\item Links of non-zero determinant whose Smith Normal Form of the coloring matrix displays different primes on different diagonal entries $($besides the $0$ entry and possible $1$'s$)$.
\item Knots whose knot group can be presented using one relator $($in particular, torus knots$)$.
\end{enumerate}
\end{cor}Proof. $1.$ and $2.$ are particular cases of Corollary \ref{cor:D}. As for $3.$, since the deficiency of knot groups is one then knot groups which can be presented with one relator only need two generators. Then the Smith Normal Form of the coloring matrix is $\text{diag } (d, 0)$ where $d$ is the determinant of the knot.
$\hfill \blacksquare$

\bigbreak

\subsection{Inner-Equivalence Classes} \label{subsect:inequivclasses}
\noindent

\begin{prop}Let $p$ be an odd prime and $n$ an integer greater than $1$.
A link $L$ with $p$-nullity $n$ has
\[
\frac{p^{n-1}-1}{2}
\]
inner-equivalence classes of $p$-colorings.
\end{prop}{\bf Proof}. As discussed right after Definition \ref{def:auto}, an inner-automorphism of $\mathbf{Z}_p$ is of the form
\[
f_{\pm, \mu} (x) = \pm  x + \mu
\]
with $\mu\in \mathbf{Z}_p$. There are then exactly $2p$ inner-automorphisms for $\mathbf{Z}_p$.

The rest of the proof goes through as in the proof of Proposition \ref{prop:numberequivclasscol} leading to the following number of inner-orbits
\[
\frac{p^n-p}{2p} = \frac{p^{n-1}-1}{2}.
\]
$\hfill \blacksquare$

\bigbreak

\begin{cor}\label{cor:innerD}
Let $L$ be a link with the following property.

The Smith Normal Form of a$($ny$)$ coloring matrix of $L$ has only one $0$ and only one $0 \neq d \neq 1$ along the diagonal.

Then for any prime $p$ such that $p|d$, there are $\frac{p-1}{2}$ equivalence class of $p$-colorings. In particular, rational  links satisfy this property.
\end{cor}Proof. Adapt the proof for Corollary \ref{cor:D}.
$\hfill \blacksquare$

\bigbreak

\begin{cor}\label{cor:innerotherD}
The following links have $\frac{p-1}{2}$ classes of $p$-colorings for each prime $p$ for which they admit non-trivial $p$-colorings.
\begin{enumerate}
\item Links whose determinant is prime.
\item Links of non-zero determinant whose Smith Normal Form of the coloring matrix displays different primes on different diagonal entries $($besides the $0$ entry and possible $1$'s$)$.
\item Knots whose knot group can be presented using one relator $($in particular, torus knots$)$.
\end{enumerate}
\end{cor}Proof. Adapt the proof for Corollary \ref{cor:otherD}.
$\hfill \blacksquare$

\bigbreak

\bigbreak

\section{Directions for Future Work} \label{sect:future}

\noindent

In the context of quandles this work has to do with homomorphisms from the fundamental quandle of the knot to the dihedral quandles (\cite{dJoyce, m}). We organize these homomorphisms into equivalence classes. In future work we plan to generalize this work to other classes of target quandles, other than the dihedral quandles.

\bigbreak

\section{Acknowledgements} \label{sect:ackn}

\noindent

J.G. thanks for support from NSFC (Grant No. 11171279 and No. 11271307).

S.J. thanks for support of the Serbian Ministry of Science (Grant No. 174012).

P.L. acknowledges support from FCT (Funda\c c\~ao para a Ci\^encia e a Tecnologia), Portugal, through project number PTDC/MAT/101503/2008, ``New Geometry and Topology''. P.L. also thanks the School of Mathematical Sciences at the University of Nottingham for hospitality.

\end{document}